\documentstyle[12pt]{article}
\newcommand\qed{{\unskip\nobreak\hfil\penalty50\hskip2em\vadjust{}
    \nobreak\hfil$\Box$\parfillskip=0pt\finalhyphendemerits=0\par}}

\renewcommand{\thesection}{\arabic{section}}
\renewcommand{\theequation}{\thesection.\arabic{equation}}

\newtheorem{theorem}{Theorem}[section]

\newtheorem{proposition}{Proposition}[section]
\newtheorem{definition}{Definition}[section]

\newcommand{\x}{\times}

\renewcommand{\a}{\alpha}
\renewcommand{\b}{\beta}
\renewcommand{\d}{\delta}
\newcommand{\D}{\Delta}
\newcommand{\e}{\varepsilon}
\newcommand{\g}{\gamma}

\newcommand{\je}{1\!\!1}

\newcommand{\pe}{p\hspace{-0.18cm}/}
\newcommand{\r}{\rho}
\newcommand{\var}{\varphi}

\newcommand{\th}{\theta}

\renewcommand{\o}{\omega}

\newcommand{\nb}{\mbox{{\bf N}}{}}
\newcommand{\rb}{\mbox{{\bf R}}{}}
\newcommand{\cb}{\mbox{{\bf C}}{}}

\newcommand{\ssb}{\mbox{{\bf S}}{}}
\newcommand{\gab}{\mbox{{\boldmath $\Gamma$}}{}}
\newcommand{\w}{\wedge}
\newcommand{\Bo}{\mbox{\raisebox{-0.3ex}{\Large\mbox{$\Box$}}}}
\newcommand{\pil}{\mbox{\large\mbox{$\pi$}}}

\newcommand{\p}{\partial}
\newcommand{\four}{(\hskip-5pt\begin{array}{c}{\scriptstyle 4}\\[-8pt]%
{\scriptstyle k}\end{array}\hskip-5pt)}

\title{Solutions of Klein--Gordon\\  and Dirac Equations on \\
Quantum Minkowski Spaces}
\author{P. Podle\'s\thanks{On leave
{}from Department of Mathematical Methods in Physics, Faculty of Physics,
University of Warsaw, Ho\.za 74, 00-682 Warszawa,
Poland}\ \thanks{This research was supported in
part by NSF grant DMS-9508597 and in part by Polish KBN grant No. 2 P301 020
07} \\
Department of Mathematics \\
University of California \\
Berkeley, CA\ \ 94720, USA}
\date{}

\begin{document}

\maketitle

\begin{abstract}
Covariant differential calculi and exterior algebras on quantum homogeneous
spaces endowed with the action of inhomogeneous quantum groups are classified.
In the case of quantum Minkowski spaces they have the same dimensions as in the
classical case.  Formal solutions of the corresponding Klein--Gordon and Dirac
equations are found.  The Fock space construction is sketched.
\end{abstract}

\setcounter{section}{-1}
\section{Introduction}
\label{sec0}

It is well known that lattice-like theories serve as regularization schemes in
quantum field theory.  But after introducing the lattice, we no longer have the
full symmetry of the original theory.  On the other hand, there was a lot of
interest in quantum spacetimes endowed with the actions of quantum groups which
are deformations of the objects used in the standard field theory (cf.
\cite{PW}, \cite{CSSW}, \cite{Lukierski}, \cite{Dobrev}, \cite{UEA},
\cite{OSWZ}, \cite{MajidP}, \cite{ChD}, \cite{Soft},
\cite{POI}).  There were two motivations of such a development:  providing
naive models of changed geometry at the Planck scale and attempts to regularize
the theory while preserving the `size' of the symmetry group in such a way that
the regularized theory could still be imagined as the theory of our universe.
Although the present paper doesn't provide support for any of these claims, we
find a lattice-like behavior of certain quantum Minkowski spaces.  It has two
aspects:

\begin{enumerate}
\item  It was found \cite{MH} that in the differential calculus
on $\rb$ corresponding
to the
one-dimensional lattice
one has
\[
xdx = (dx)x + l dx
\]
where $x$ is the identity function and $l$
is the lattice constant.  In
Section~\ref{sec1} we describe
 differential calculi on quantum Minkowski spaces
by a very similar relation (\ref{eq1.3}).

\item  For the above differential calculus on $\rb$ one has
$df=dx\p(f)=\tilde\p(f)dx$, where $f$ is a function on $\rb$,
$\p(f)=(f(x+l)-f(x))/l$,
$\tilde\p(f)=(f(x)-f(x-l))/l$ (cf. \cite{MH}).
Setting $\D=\p\tilde\p=\tilde\p\p$, one gets
\[
\D e^{-ipx} = \frac{\mbox{sin}^2k}{k^2}(-p^2)e^{-ipx},
\]
where $k=pl/2$. Thus we obtain
an additional factor
$\frac {\sin^2 k}{k^2}$
(comparing with the action of the usual Laplacian $\p^2/\p x^2$).
In Section~\ref{sec4} similar factors appear in the
description of eigenvalues of
the Laplacian on quantum Minkowski spaces.
\end{enumerate}

In Section~\ref{sec1} we recall the definition of homogeneous quantum spaces
$M$ (e.g. quantum Minkowski spaces) endowed with the action of inhomogeneous
quantum groups $G$ (e.g. quantum Poincar\'e groups).  We classify the
differential calculi on $M$ which have the same properties as in the classical
case.  They exist if and only if a certain matrix ${\tilde F}$
(related to the existence of quasitriangular structure on $G$ \cite{CQT})
vanishes,
in which case they are unique.
  In Section~\ref{sec2} we prove that each such calculus has a unique
natural extension to an exterior algebra of differential forms.  In the case of
quantum Minkowski spaces the modules of $k$-forms have the classical dimensions
$\four$.  In Section~\ref{sec3}
the properties of partial derivatives, Laplacian and the Dirac operator are
investigated.  We heuristically assert hermiticity
of the momenta and Laplacian.  In Section~\ref{sec4}
we find formal solutions of Klein--Gordon and Dirac equations for two special
classes of $M$.  They are obtained from the plane waves $e^{-ip_ax^a}$, but
the
eigenvalues of the
momenta are related to $p_a$ in a complicated
way in general.  The sketch of the Fock space construction is provided in
Section~\ref{sec5}.

We sum over repeated indices (Einstein's convention).  We use covariant and
contravariant indices but each object is used only with one fixed position of
indices (except for $g^{ab}$ and $g_{ab}$
which are inverse one to another:  $g^{ab}g_{bc} = g_{cb}g^{ba} = {\d^a}_c$).
  If
$V,W$ are vector spaces then $\tau: V \otimes W \to W \otimes V$ is given by
$\tau(x \otimes y) = y \otimes x$, $x \in V$, $y \in W$.  We denote the unit
matrix by $1\!\!1$, $1\!\!1^{\otimes k} = 1\!\!1 \otimes \ldots \otimes 1\!\!1$
($k$ times).  If ${\cal A}$ is an algebra, $v \in M_N({\cal A})$, $w \in
M_K({\cal A})$, then the tensor product $v \otimes w \in M_{NK}({\cal A})$ is
defined by
\[
(v \otimes w)^{ij}{}_{kl} = v^i{}_kw^j{}_l,\ i,k = 1,\dots,N,\ %
j,l = 1,\dots,K.
\]
We set $\dim v = N$.  If ${\cal A}$ is a ${}^*$-algebra then
the conjugate of $v$ is defined as ${\bar v} \in M_N({\cal A})$ where ${\bar
v}^i{}_j = (v^i{}_j)^*$.  We also set $v^* = {\bar v}^T$ ($v^T$ denotes the
transpose of $v$, i.e. $(v^T)_i{}^j = v^j{}_i$).

Throughout the paper quantum groups $H$ are abstract objects described by the
corresponding Hopf (${}^*$-) algebras $\mbox{Poly}(H) = ({\cal A},\D)$.  We
denote
by $\D,\e,S$ the comultiplication, counit and the coinverse of
$\mbox{Poly}(H)$.  We say that $v$ is a representation of $H$ (i.e. $v \in
\mbox{ Rep } H$) if $v \in M_N({\cal A})$, $N \in {\nb}$, and
\[
\D v^i{}_j = v^i{}_k \otimes v^k{}_j,\ \e(v^i{}_j) =
{\d^i}_j,\ i,j = 1,\dots,N,
\]
in which case
 $S(v^i{}_j)=(v^{-1})^i{}_j$. Matrix elements of all $v\in\mbox{ Rep } H$
linearly span ${\cal A}$.
The conjugate of a representation and tensor products of representations are
also representations.  The set of nonequivalent irreducible representations of
$H$ is denoted by $\mbox{Irr } H$.  If $v,w \in \mbox{ Rep } H$, then we say
that $A \in M_{\dim v \x \dim w}({\cb})$ intertwines $v$ with $w$ (i.e. $A
\in \mbox{ Mor}(v,w)$) if $Av = wA$.  For $\rho,\rho' \in {\cal A}'$ (the dual
vector space of ${\cal A}$) one defines their convolution $\rho * \rho' = (\rho
\otimes \rho')\D$. For $\rho\in {\cal A}'$, $a\in {\cal A}$, we set
$\rho * a=( id \otimes\rho)\D a$, $a*\rho=(\rho\otimes id )\D a$.

\section{The covariant differential calculi on quantum homogeneous spaces}
\label{sec1}

In this section we recall the definition of the quantum homogeneous space $M$
endowed with the action of inhomogeneous quantum group $G$ \cite{INH}.  The
corresponding unital algebra ${\cal C} = \mbox{ Poly}(M)$ is generated by
quantum coordinates $x^i$, $i = 1,\dots,N$.  We prove that there exists a
covariant differential calculus on $M$ which has $dx^i$, $i = 1,\dots,N$, as
the basis of the module of $1$-forms if and only if a certain matrix ${\tilde
F} = 0$.  Moreover, such a calculus is unique.  We specify
the quantum Minkowski
spaces endowed with the action of quantum Poincar\'e groups \cite{POI} for
which ${\tilde F} = 0$.

Throughout the section $\mbox{Poly}(H) = ({\cal A},\D)$ is any Hopf algebra
with invertible $S$ (we need invertibility of $S$ in the proof of
Theorem \ref{th1.2}, it was not needed in \cite{INH}) such that

\begin{enumerate}
\begin{enumerate}
\item  each representation of $H$ is completely reducible,
\item  $\Lambda$ is an irreducible representation of $H$,

\item  $\mbox{Mor}(v \otimes w,\ \Lambda \otimes v \otimes w) = \{0\}$ for any
two irreducible representations of $H$.
\end{enumerate}
\end{enumerate}

Moreover, we assume that $f^i{}_j$, $\eta^i \in {\cal A}'$, $T^{ij} \in {\cb}$,
 $i,j = 1,\dots,N = \dim \Lambda$, are given and satisfy

\begin{enumerate}
\item ${\cal A} \ni a \to \rho(a) = \left( \begin{array}{cc} f(a) & \eta(a) \\
0 & \e(a) \end{array} \right) \in M_{N+1}({\cb})$ is a unital homomorphism,

\item $\Lambda^s{}_t(f^t{}_r * a) = (a * f^s{}_t)\Lambda^t{}_r$ for $a \in
{\cal A}$,

\item $R^2 = 1\!\!1$ where $R^{ij}{}_{sm} = f^i{}_m(\Lambda^j{}_s)$,

\item $(\Lambda \otimes \Lambda)^{kl}{}_{ij}(\tau^{ij} * a) = a * \tau^{kl}$
for $a \in {\cal A}$ where $\tau^{ij} = (R - 1\!\!1)^{ij}{}_{mn}(\eta^n *
\eta^m - \eta^m(\Lambda^n{}_s)\eta^s + T^{mn} \e -f^n{}_b * f^m{}_aT^{ab})$,

\item $A_3{\tilde F} = 0$ where
$A_3 = 1\!\!1 \otimes 1\!\!1 \otimes 1\!\!1 - R \otimes 1\!\!1 - 1\!\!1 \otimes
R + (R \otimes 1\!\!1)(1\!\!1 \otimes R) + (1\!\!1 \otimes R)(R \otimes 1\!\!1)
- (R \otimes 1\!\!1)(1\!\!1 \otimes R)(R \otimes 1\!\!1)$, ${\tilde
F}^{ijk}{}_m = \tau^{ij}(\Lambda^k{}_m)$,

\item $A_3(Z \otimes 1\!\!1 - 1\!\!1 \otimes Z)T = 0$, $RT = -T$ where
$Z^{ij}{}_k = \eta^i(\Lambda^j{}_k)$.
\end{enumerate}

In particular, 4.--5. are satisfied if $\tau^{ij} = 0$.  The inhomogeneous
quantum group $G$ corresponds to the Hopf algebra $\mbox{Poly}(G) = ({\cal
B},\D)$ defined (cf. Corollary $3.8.a$ of \cite{INH}) as follows:  ${\cal
B}$ is the universal unital algebra generated by ${\cal A}$ and $y^i$, $i =
1,\dots,d$, satisfying the relations $I_{\cal B} = I_{\cal A}$,
\setcounter{equation}{0}
\renewcommand{\theequation}{\thesection.\arabic{equation}}
\begin{equation}
\label{eqA}
y^sa = (a *
f^s{}_t)y^t + a * \eta^s -\Lambda^s{}_t(\eta^t * a),\ a \in {\cal A},
\end{equation}
\begin{equation}
\label{eqB}
(R
- 1\!\!1)^{kl}{}_{ij}(y^iy^j - \eta^i(\Lambda^j{}_s)y^s + T^{ij} -
\Lambda^i{}_m\Lambda^j{}_nT^{mn}) = 0.
\end{equation}
Moreover, $({\cal A},\D)$ is a Hopf subalgebra of $({\cal B},\D)$ and $\D y^i =
\Lambda^i{}_j \otimes y^j + y^i \otimes I$
(note that the $y^i$ were denoted by $p_i$ in
\cite{INH}).  We define
$\cal C = \mbox{Poly}(M)$ as the universal unital algebra
generated by $x^i$, $i = 1,\dots,N$, satisfying
\begin{equation}
\label{eq1.1}
(R - 1\!\!1)^{ij}{}_{kl}(x^kx^l - Z^{kl}{}_s x^s + T^{kl}) = 0.
\end{equation}
The action of $G$ on $M$ is
described by the unital homomorphism $\Psi: {\cal C} \to {\cal B} \otimes {\cal
C}$ such that
\begin{equation}
\label{eq1.2}
\Psi(x^i) = \Lambda^i{}_j \otimes x^j + y^i \otimes I,
\end{equation}
$(\e \otimes id)\Psi = id,\ (id \otimes \Psi)\Psi = (\D \otimes id)\Psi.$
The pair $({\cal C},\Psi)$ was investigated in Section~5 of \cite{INH}.  We
assume
\begin{equation}
\label{eq1.0}
\mbox{Mor}(I,\Lambda \otimes \Lambda \otimes \Lambda) = \{0\},
\end{equation}
\begin{equation}
\label{eq1.0'}
\mbox{Mor}(I,\Lambda \otimes \Lambda) \cap \ker(R + 1\!\!1) = \{0\}
\end{equation}
((\ref{eq1.0})--(\ref{eq1.0'}) are satisfied for $G$ being quantum Poincar\'e
groups \cite{POI}).
Then $M$ is
called quantum homogeneous space and has the properties analogous to the
Minkowski space (cf. Section~5 of \cite{INH}, Section~1 of \cite{POI}).
The `sizes' of ${\cal B}$ and ${\cal C}$ were described in Corollary 3.6
and Proposition 5.3 of \cite{INH}.

Motivated by \cite{DC}, \cite{TSQ}, \cite{DCQS} we have

\begin{definition}
\label{def1.1}
We say that ${{\gab}}^{\w 1} = (\Gamma^{\w 1},\Psi^{\w
1},d)$ is a covariant differential calculus on $M$ if

\begin{enumerate}
\item $\Gamma^{\w 1}$ is a ${\cal C}$-bimodule,
$\o I_{\cal C} = I_{\cal C}\o
= \o$ for $\o \in \Gamma^{\w 1}$,

\item $\Psi^{\w 1}: \Gamma^{\w 1} \to {\cal B} \otimes
\Gamma^{\w 1}$ satisfies

\begin{enumerate}
\item $(\e \otimes id)\Psi^{\w 1} = id$, $(id \otimes \Psi^{\w
1})\Psi^{\w 1} = (\D \otimes id)\Psi^{\w 1}$,

\item $\Psi^{\w 1}(\o a) = \Psi^{\w 1}(\o)\Psi(a)$, $\Psi^{\w
1}(a\o) = \Psi(a)\Psi^{\w 1}(\o)$ for $\o \in \Gamma^{\w 1}$, $a \in
{\cal C}$,
\end{enumerate}

\item $d: {\cal C} \to \Gamma^{\w 1}$ is a linear map such that

\begin{enumerate}
\item $d(ab) = a(db) + (da)b$, $a,b \in {\cal C}$,

\item $(id \otimes d)\Psi = \Psi^{\w 1}d$,

\item $\Gamma^{\w 1} = \mbox{ {\em span}}\{(da)b: a,b \in {\cal C}\}$.
\end{enumerate}
\end{enumerate}

We say that $\Gamma^{\w 1}$ is $N$-dimensional if $dx_i$, $i = 1,\dots,N$,
form a basis of $\Gamma^{\w 1}$ (as right ${\cal C}$-module).
\end{definition}

\begin{theorem}
\label{th1.2}
There exists $N$-dimensional covariant differential calculus on $M$ iff
${\tilde F} = 0$.  In that case it is uniquely determined by
\begin{equation}
\label{eq1.3}
x^idx^j = {R^{ij}}_{kl} dx^kx^l + Z^{ij}{}_kdx^k,\ i,j = 1,\dots,N,
\end{equation}
\begin{equation}
\label{eq1.4}
\Psi^{\w 1}dx^i = \Lambda^i{}_j \otimes dx^j,\ i = 1,\dots,N.
\end{equation}
\end{theorem}

\medskip
\noindent
{\bf Proof}.  Let ${{\gab}}^{\w 1}$ be $N$-dimensional covariant differential
calculus on $M$.  Using (\ref{eq1.2}) and condition 3b) of
Definition~\ref{def1.1}, one gets (\ref{eq1.4}).  The linear mappings $\p_i:
{\cal C} \to {\cal C}$, $\r_i{}^j: {\cal C} \to {\cal C}$, $i,j = 1,\dots,N$,
are uniquely defined by
\begin{equation}
\label{eq1.5}
da = dx^i\p_i(a),\ adx^i = dx^j\r_j{}^i(a).
\end{equation}
Using condition 3a) of Definition~\ref{def1.1} and $(ab)dx^i = a(bdx^i)$,
one gets that
\begin{equation}
\label{eq1.5a}
{\cal L}: {\cal C} \ni a \to \left[ \begin{array}{cc}
(\r_i{}^j(a))_{i,j=1}^N & (\p_i(a))_{i=1}^N \\
0 & a
\end{array} \right] \in M_{N+1}({\cal C})
\end{equation}
is a unital homomorphism.  Conditions
 2b) and 3b) of Definition~\ref{def1.1}
imply
\begin{equation}
\label{eq1.6}
(id \otimes \p_j)\Psi(a) = (\Lambda^i{}_j \otimes I)\Psi(\p_i(a)),\ j =
1,\dots,N,
\end{equation}
\begin{equation}
\label{eq1.7}
(id \otimes \r_k{}^j)\Psi(a)(\Lambda^i{}_j \otimes I) = (\Lambda^j{}_k \otimes
I)\Psi(\r_j{}^i(a)),\ i,k = 1,\dots,N,
\end{equation}
$a \in {\cal C}$.  Moreover,
\begin{equation}
\label{eq1.8}
\p_j(x^i) = \d^i{}_j,\ i,j = 1,\dots,N.
\end{equation}
Conversely, any unital homomorphism (\ref{eq1.5a}) satisfying
(\ref{eq1.6})--(\ref{eq1.8}) for $a \in {\cal C}$ defines,
 through (\ref{eq1.5}),
a covariant $N$-dimensional differential calculus on $M$.  So we need to find
all such homomorphisms.  Let us notice that it is sufficient to check
(\ref{eq1.6})--(\ref{eq1.7}) for $a = x^l$, $l = 1,\dots,N$ (they are trivial
for $a = I$ and if $a,b$ satisfy them then---using the
homomorphism property---$a
\cdot b$ does
also).  But for $a = x^l$ (\ref{eq1.6})
is trivial.  Moreover, ${\cal
L}$ is determined by
\[
{\cal L}(x^l) = h^l \equiv \left[ \begin{array}{cc}
(K_i{}^{lj})_{i,j=1}^N & (\d^l{}_i)_{i=1}^N \\
0 & x^l \end{array} \right],
\]
$l = 1,\dots,N$, where $K_i{}^{lj} = \r_i{}^j(x^l)$.  The equation
(\ref{eq1.8}) and the existence of ${\cal L}$ (for given $K$) are equivalent
to (\ref{eq1.1}) with $x^i$ replaced by $h^i$.  This can be translated to
\begin{equation}
\label{eq1.9}
d[(R - 1\!\!1)^{ij}{}_{kl}(x^kx^l - Z^{kl}{}_s x^s + T^{kl})] = 0,
\end{equation}
\begin{equation}
\label{eq1.10}
[(R - 1\!\!1)^{ij}{}_{kl}(x^kx^l - Z^{kl}{}_s x^s + T^{kl})]dx^m = 0,
\end{equation}
where the left hand sides should be expanded using condition 3a) of
Definition~\ref{def1.1} and (\ref{eq1.5}) so that $dx^s$ appear on the very
left of the equations. Then the condition means that the total coefficient
multiplying $dx^s$ from the right is zero.  The condition (\ref{eq1.7})
for $a = x^l$ means
\[
\Lambda^l{}_t\Lambda^i{}_j \otimes K_{k}{}^{tj} + y^l\Lambda^i{}_k \otimes I =
(\Lambda^j{}_k \otimes I)\Psi(K_j{}^{li}),
\]
i.e.
\begin{equation}
\label{eq1.11}
\Psi(K_m{}^{li}) = (G \otimes \Lambda \otimes \Lambda)_m{}^{li},{}^k{}_{tj}
\otimes K_k{}^{tj} + G_m{}^ky^l\Lambda^i{}_k \otimes I
\end{equation}
where
$(G^{-1})_i{}^j = \Lambda^j{}_i$ ($G = [S^{-1}(\Lambda)]^T =
(\Lambda^T)^{-1}$).
Thus we need to find $K$ satisfying (\ref{eq1.9})--(\ref{eq1.11}).

Now (\ref{eq1.3}) is equivalent to $K_k{}^{ij} = \r_k{}^j(x^i) = R^{ij}{}_{kl}
x^l + Z^{ij}{}_k$.  It is easy to check that such a
$K$ satisfies (\ref{eq1.9})
and (\ref{eq1.11}).  Suppose there exists another ${\tilde K}$ satisfying
(\ref{eq1.11}).  Then $M = K - {\tilde K}$ satisfies $\Psi(M_m{}^{li}) = (G
\otimes \Lambda \otimes \Lambda)_m{}^{li},{}^k{}_{tj} \otimes M_k{}^{tj}$.
Using Condition~2. of Section~5 of \cite{INH}, one gets $M_m{}^{li} \in
{\cb}$,  $G_m{}^k\Lambda^l{}_t\Lambda^i{}_jM_k{}^{tj} = M_m{}^{li}$.
Multiplying from the left by $\Lambda^m{}_s = (G^{-1})_s{}^m$ and setting
$U^{tj}{}_s = M_s{}^{tj}$, one gets $\Lambda^l{}_t\Lambda^i{}_jU^{tj}{}_s =
U^{li}{}_m\Lambda^m{}_s$, i.e.
\[
U \in \mbox{ Mor}(\Lambda,\Lambda \otimes \Lambda) = \{0\},\ U = 0,\ M = 0,\
{\tilde K} = K.
\]
Therefore uniqueness follows.  Expanding (\ref{eq1.10}) and using
(\ref{eq1.3}), one gets that (after a long computation) the total coefficients
multiplying $dx^s$ from the right are zero if and only if ${\tilde F} = 0$
(we use the results of \cite{INH}: Proposition 5.3 and Remark 5.4 for
$N=1$, (3.61) and (3.30)).
Thus the existence statement and (\ref{eq1.3}) are proved. \qed

\medskip
All the assumptions (including (\ref{eq1.0})--(\ref{eq1.0'})) are fulfilled if
$H$ is a quantum Lorentz group \cite{WZ}, $G$ is a quantum Poincar\'e group
and
$M$ is the corresponding unique quantum Minkowski space \cite{POI}.  According
to Theorem~\ref{th1.2}, there exists 4-dimensional covariant differential
calculus on the quantum Minkowski space iff ${\tilde F} = 0$,
 iff $\lambda = 0$
(see the proof of Theorem~$1.6$ of \cite{POI}),
 which holds for all cases except
of the following:

\noindent
1), $t = 1$, $s = 1$, $t_0 \in {\rb}{\backslash}\{0\}$

\noindent
5), $t = 1$, $s = \pm 1$, $t_0 \in {\rb}{\backslash}\{0\}$

\noindent
4), $s = 1$, $b \ne 0$

\noindent
(in the terminology of Theorem~$1.6$ and Remark~$1.8$ of \cite{POI}\footnote{
In the old version of Remark 1.8 of \cite{POI} one should replace $t$
by $t_0$ (except of expressions $t=1$). That $t_0$ is identified with
$t$ of (3)--(4) of Ref. 16 of \cite{POI}.}).  Then $N = 4$.
Such a calculus is unique.

Let $({\cal A},\D)$ be a Hopf ${}^*$-algebra.  Then $S$ is always invertible.
We also assume
\[
{\bar \Lambda} = \Lambda,\ f^i{}_j(S(a^*)) = \overline{f^i{}_j(a)},\
\eta^i(S(a^*)) = \overline{\eta^i(a)},\ \overline{T^{ij}} = T^{ji},\ i,j =
1,\dots,N,\ a \in {\cal A}.
\]
In that case \cite{INH} $({\cal B},\D)$ has a unique Hopf ${}^*$-algebra
structure
such that $({\cal A},\D)$ is its Hopf ${}^*$-subalgebra and $y^{i*} = y^i$.
Moreover, ${\cal C}$ is a ${}^*$-algebra with $*$ defined by $x^{i*} = x^i$ and
$\Psi$ is a ${}^*$-homomorphism.

\begin{proposition}
\label{prop1.3}
Under assumptions as above the $N$-dimensional covariant differential calculus
on $M$ described in Theorem~\ref{th1.2} possesses a unique $*: \Gamma^{\w
1} \to \Gamma^{\w 1}$ such that:

\begin{enumerate}
\item $(\o a)^* = a^*{\o}^*$, $(a\o)^* = {\o}^*a^*$,
$\o \in \Gamma^{\w 1}$, $a \in
{\cal C}$,

\item $(da)^* = d(a^*)$, $a \in {\cal C}$,

\item $\Psi^{\w 1}({\o}^*) = (\Psi^{\w 1}(\o))^{*\otimes *}$, $\o \in
\Gamma^{\w 1}$.
\end{enumerate}
\end{proposition}

\medskip
\noindent
{\bf Proof}.  We must define $*: \Gamma^{\w 1} \to \Gamma^{\w 1}$ by
\begin{equation}
\label{eq1.12}
(dx^ia_i)^* = a_i{}^*dx^i,\ a_i \in {\cal C}.
\end{equation}
In virtue of $(\ref{eq1.3})$ of the present paper,
 $(4.14)$ and the next formula of \cite{INH}
\[
\begin{array}{rll}
(x^idx^j)^* &= &(R^{ij}{}_{kl}dx^kx^l + Z^{ij}{}_k dx^k)^* \\
&= &R^{ji}{}_{lk}(x^ldx^k - Z^{lk}{}_sdx^s) = dx^jx^i.
\end{array}
\]
But
\[
\begin{array}{rll}
(x^idx^j)^* &= &(dx^k\r_k{}^j(x^i))^* \\
&= &\r_k{}^j(x^i)^*dx^k = dx^s[\r_s{}^k(\r_k{}^j(x^i)^*)].
\end{array}
\]
Therefore $\r_s{}^k(\r_k{}^j(a)^*) = a^*{\d^j}_s$
 for $a = x^i$ and hence ($\r$ is a
unital homomorphism) for all $a \in {\cal A}$.  This means $(a_jdx^j)^* =
dx^ja_j^*$, $a_j \in {\cal C}$.  This and (\ref{eq1.12}) prove
condition~$1.$ for any $\o = dx^ia_i \in \Gamma^{\w 1}$.  Writing $a \in
{\cal C}$ as a polynomial in $x^j$, using
condition~$1.$ and
 $(dx^i)^* = dx^i$ (see
(\ref{eq1.12})),
one gets condition~$2.$.  By virtue
of (\ref{eq1.12}) and (\ref{eq1.4}) we obtain condition~$3.$. \qed

\medskip
In particular, all the above ${}^*$-structures exist for
quantum Poincar\'e groups \cite{POI}, quantum Minkowski spaces
\cite{POI} and $4$-dimensional covariant differential calculi on them.

\medskip
{\bf Remark.}
In the case of $Z=T=0$ formulae (\ref{eq1.3}), (\ref{eq2.2}), (\ref{eq3.1}),
(\ref{eq3.1''}),
(\ref{eq3.2}), (\ref{eq3.6}) and the second formula of (\ref{eq3.7}) or their
analogues were studied in several contexts in \cite{TSQ},
\cite{WeZu}, \cite{CSW}, \cite{OSWZ}, \cite{UEA}, \cite{ChD}, \cite{Kulish}.

\section{Exterior algebras}
\label{sec2}

In this section we construct the exterior algebras for the $N$-dimensional
covariant differential calculi described in the previous section.  In the case
of quantum Minkowski spaces the right ${\cal C}$-module of $k$-forms has
dimension
$\four$ as in the
classical case.

Throughout this section ${{\gab}}^{\w 1} = (\Gamma^{\w
1},\Psi^{\w 1},d)$ is an $N$-dimensional covariant differential calculus
on quantum homogeneous space $M$ endowed with the action of inhomogeneous
quantum group $G$ as described in Theorem~\ref{th1.2}.  In particular, we
assume
all the conditions introduced before
Theorem~\ref{th1.2} and that
${\tilde F} = 0$.

\begin{definition}[cf. \cite{DC}, \cite{TSQ}, \cite{DCQS}]
\label{def2.1}
We say that ${{\gab}}^{\w} = (\Gamma^{\w},\Psi^{\w},d)$ is
an exterior algebra on $M$ iff

\begin{enumerate}
\item $\Gamma^{\w} = \oplus_{n=0}^{\infty} \Gamma^{\w n}$ is a graded
algebra such that $\Gamma^{\w 0} = {\cal C}$ and the unit of ${\cal C}$
is the unit of $\Gamma^{\w}$,

\item $\Psi^{\w}: \Gamma^{\w} \to {\cal B} \otimes \Gamma^{\w}$
is a graded homomorphism such that
\[
(\e \otimes id)\Psi^{\w} = id,\quad (id
\otimes \Psi^{\w})\Psi^{\w} = (\D \otimes id)\Psi^{\w},\quad
\Psi^{\w 0} = \Psi,
\]

\item $d: \Gamma^{\w} \to \Gamma^{\w}$ is a linear mapping such that

\begin{enumerate}
\item $d(\Gamma^{\w n}) \subset \Gamma^{\w(n+1)}$, $n = 0,1,2,\dots$,

\item $d(\th \wedge \th') = d\th \wedge \th' + (-1)^k\th \wedge d\th'$, $\th
\in \Gamma^{\w k}$, $\th' \in \Gamma^{\w}$ ($\wedge$ denotes
multiplication in $\Gamma^{\w}$),

\item $(id \otimes d)\Psi^{\w} = \Psi^{\w}d$,

\item $dd = 0$,
\end{enumerate}

\item $\Gamma^{\w n} = \mbox{ span}\{(da_1 \wedge \ldots \wedge da_n)a_0:
a_0,a_1,\dots,\ a_n \in {\cal C}\}$ (we omit $\wedge$ if one of multipliers
belongs to ${\cal C}$),

\item $\Gamma^{\w 1}$, $\Psi^{\w 1}$, $d: {\cal C} \to
\Gamma^{\w 1}$ are as in Definition~\ref{def1.1} and Theorem~\ref{th1.2},

\item if $({\tilde \Gamma}^{\w},{\tilde \Psi}^{\w},{\tilde d})$ also
satisfies 1.--5. then there exists a graded homomorphism $\rho:
\Gamma^{\w} \to {\tilde \Gamma}^{\w}$ which is an identity on ${\cal
C}$ and satisfies ${\tilde \Psi}^{\w}\rho = (id \otimes
\rho)\Psi^{\w}$, ${\tilde d}\rho = \rho d$ (universality condition).
\end{enumerate}
\end{definition}

We set $R_{nk} = 1\!\!1^{\otimes(k-1)} \otimes R \otimes
1\!\!1^{\otimes(n-k-1)}$, $R_{n\pi} = R_{nk_1} \cdot \dots \cdot R_{nk_s}$ for
any permutation $\pi = t_{k_1} \cdot \dots \cdot t_{k_s} \in \Pi_n$ where $t_k$
is the transposition $k \leftrightarrow k+1$, $A_n = \frac {1}{n!} \sum_{\pi
\in \Pi_n} (-1)^{\mbox{sgn } {\pil}}
R_{n\pi}$, $A_n^2 = A_n$, $R_{nk}A_n = A_nR_{nk}
= -A_{n}$, $k=1,2,\dots,n-1$.

Let ${\a}' = \{\a'_i: i = 1,\dots,\dim A_n\}$ be a basis of $\mbox{im }
A_n$, $\b' = \{\b'_j: j = 1,\dots,\dim(1\!\!1 - A_n)\}$ be a basis of
$\mbox{im}(1\!\!1 - A_n)$.  Then $\a' \sqcup \b'$ is a basis of
$({\cb}^N)^{\otimes n}$.  We denote by $\a \sqcup \b$ the dual basis.
Therefore
\setcounter{equation}{0}
\begin{equation}
\label{eq2.0}
\a^i A_n = \a^i,\ \b^j  A_n = 0,
\end{equation}
\begin{equation}
\label{eq2.0'}
A_n = \a'_i\a^{i}.
\end{equation}

\begin{theorem}
\label{th2.2}
There exists a unique exterior algebra ${{\gab}}^{\w}$ on $M$.  The
$n$-forms
\begin{equation}
\label{eq2.1}
{\o}^{\g} = {\a}_{k_1,\dots,k_n}^{\g n} dx^{k_1} \wedge \ldots \wedge
dx^{k_n},\ \g = 1,\dots,\dim A_n,
\end{equation}
form a basis of the right ${\cal C}$-module $\Gamma^{\w n}$.  Moreover,
\begin{equation}
\label{eq2.2}
dx^i \wedge dx^j = -R^{ij}{}_{kl}dx^k \wedge dx^l,\ i,j = 1,\dots,N.
\end{equation}
\end{theorem}

\medskip
\noindent
{\bf Proof}.  Assume that ${\ssb}^{\w} =
(S^{\w},\Psi_S^{\w},d_S)$ satisfies conditions 1.--5. of
Definition~\ref{def2.1}.  Acting $d_S$ on (\ref{eq1.3}), one gets
(\ref{eq2.2}).  Set $dx_S^J = dx^{i_1} \wedge \ldots \wedge dx^{i_n}$ where $J
= (i_1,\dots,i_n) \in \{1,\dots,N\}^n \equiv N^{(n)}$.  They generate the right
${\cal C}$-module $\Gamma^{\w n}$.  We obtain $dx_S^J =
-(R_{nk})^J{}_Kdx_S^K$, $dx_S^J = (-1)^{\mbox{sgn }
{\pil}}
(R_{n\pi})^J{}_Kdx_S^K$, $dx_S^J = (A_n)^J{}_Kdx_S^K$, (\ref{eq2.1})
generate the right ${\cal C}$-module $S^{\w n}$.  Moreover, $S^{\w n}
= \{dx_S^Ja_J: J \in N^{(n)}\}$, formulae for $d_S(dx_S^Ja_J)$,
$\Psi_S^{\w}(dx_S^Ja_J)$ are determined by conditions 2. and
3(b)(c)(d)
of Definition~\ref{def2.1}.  Thus it suffices to construct
${{\gab}}^{\w} = (\Gamma^{\w},\Psi^{\w},d)$ which
satisfies the conditions 1.--5. of Definition~\ref{def2.1} and such that
${\o}^{\g}$, $\g = 1,\dots,\dim A_n$, are independent in the right ${\cal
C}$-module $\Gamma^{\w n}$.

We set $\Gamma = \Gamma^{\w 1}$, $\Gamma^{\otimes 0} = {\cal C}$,
$\Gamma^{\otimes n} = \Gamma \otimes_{\cal C} \ldots \otimes_{\cal C} \Gamma$,
$\Gamma^{\otimes} = \oplus_{n=0}^{\infty} \Gamma^{\otimes n}$, $dx^I =
dx^{i_1} \wedge \ldots \wedge dx^{i_n}$ for $I = (i_1,\dots,i_n) \in
N^{(n)}$,
\[
L_0^n = \mbox{ span}\{dx^J - (A_n)^J{}_K dx^K: J \in N^{(n)}\},
\]
\[
L^n = L_0^n{\cal C},\ L = \oplus_{n=0}^{\infty} L^n,\ \Gamma^{\w n} =
\Gamma^{\otimes n}/L^n,\ \Gamma^{\w} = \oplus_{n=0}^{\infty}
\Gamma^{\w n} = \Gamma^{\otimes}/L.
\]
Then $\Gamma^{\otimes n}$ are ${\cal C}$-bimodules, $\Gamma^{\otimes}$ is a
graded algebra, ${\o}^{\g}$, $\g = 1,\dots,\dim A_n$, form a basis of the right
${\cal C}$-module $\Gamma^{\w n}$ (cf. (\ref{eq2.0})).  We see that
$L_0^n = \mbox{ span}\{dx^J + (R_{nk})^J{}_Kdx^K: J \in N^{(n)},\ k =
1,\dots,N-1\}$.  Thus $L$ is the right ideal generated by $dx^K \otimes
\o^{ik}$, $K \in N^{(s)}$, $s = 0,1,\dots$, $i,k = 1,\dots,N$, where $\o^{ik} =
dx^i \otimes_{\cal C} dx^j + R^{ij}{}_{kl}dx^k \otimes_{\cal C} dx^l$.  Using
(\ref{eq1.3}) and $(3.11)$, $(3.61)$ of \cite{INH}, one shows that
\[
x^m\o^{ik} = R^{sk}{}_{nb} R^{mi}{}_{js} \o^{jn}x^b + \o^{ab}l^{mik}{}_{ab},
\]
where $l = Z \otimes 1\!\!1 + (R \otimes 1\!\!1)(1\!\!1 \otimes Z)$.  That and
(\ref{eq1.3}) prove that $L$ is the ideal generated by $\o^{ik}$ and
condition 1. of Definition~\ref{def2.1} follows.  Moreover, (\ref{eq2.2}) is
satisfied.

We define a linear mapping $\Psi^{\otimes}: \Gamma^{\otimes} \to
\Gamma^{\otimes}$ by
\begin{equation}
\label{eq2.2'}
\Psi^{\otimes}(\o_1 \otimes_{\cal C} \ldots \otimes_{\cal
C} \o_n) = \o_1^{(1)}\o_2^{(1)} \cdot \dots \cdot \o_n^{(1)} \otimes \o_1^{(2)}
\otimes_{\cal C} \ldots \otimes_{\cal C} \o_n^{(2)},
\end{equation}
where $\o_s \in \Gamma^{\w 1}$, $\Psi^{\w 1}(\o_s) = \o_s^{(1)} \otimes
\o_s^{(2)}$ (Sweedler's notation), $s = 1,\dots,n$.  Then $\Psi^{\otimes}$ is
well defined, $\Psi^{\otimes} \o^{ik} = \Lambda^i{}_j\Lambda^k{}_m \otimes
\o^{jm}$ (see (\ref{eq1.4})), $\Psi^{\otimes}(L) \subset {\cal B} \otimes L$,
$\Psi^{\otimes}$ defines $\Psi^{\w}: \Gamma^{\w} = \Gamma^{\otimes}/L
\to {\cal B}\otimes\Gamma^{\w}$ which satisfies condition 2. of
Definition~\ref{def2.1}.

We set $d(\o^{\g}a_{\g})=(-1)^n\o^{\g}\w da_{\g}$, i.e. (see (\ref{eq2.1}),
(\ref{eq2.0'}))
\begin{equation}
\label{eq2.3}
d(dx^Ja_J) = (-1)^{|J|}dx^J \wedge da_J,\ a_J \in {\cal C},\ J \in N^{(n)},\
|J| = n.
\end{equation}
Conditions 5., 4. and 3(a)(c) follow.  Moreover,
\begin{equation}
\label{eq2.3'}
ddx^J = 0.
\end{equation}
We shall prove
\begin{equation}
\label{eq2.4}
d(a_Jdx^J) = da_J \wedge dx^J.
\end{equation}
Due to (\ref{eq1.3}), (\ref{eq2.3}) and (\ref{eq2.2})
\[
\begin{array}{rll}
d(x^idx^j) &= &d[R^{ij}{}_{kl}(dx^k)x^l + Z^{ij}{}_ldx^l] \\
&= &-R^{ij}{}_{kl}dx^k \wedge dx^l = dx^i \wedge dx^j.
\end{array}
\]
Thus $x^i \in E = \{a \in {\cal C}: d(adx^i) = da \wedge dx^i$, $i =
1,\dots,N\}$.  Moreover, $E$ is a unital algebra (we use (\ref{eq1.5}) and
(\ref{eq2.3}) in order to find LHS in the definition of $E$ and perform
a direct computation).
Hence (\ref{eq2.4}) for $|J| = 1$ follows.
According to (\ref{eq2.3}), $d(dx^i\w\o)=-dx^i\w d\o$ for $\o\in\Gamma^{\w}$.
 Using this, (\ref{eq1.5}) and the mathematical induction
w.r.t. $|J|$, a simple calculation proves (\ref{eq2.4}) in the general case.
Then it is easy to check condition 3(b).

We set $F = \{a \in {\cal C}: dda = 0\}$.  Then $x^i \in F$ (see
(\ref{eq2.3'})) and $F$ is a unital algebra (we use (\ref{eq2.3}),
(\ref{eq2.4})).  Thus $F = {\cal C}$.  This and (\ref{eq2.3'}) show $dd\o = 0$
for any $\o = dx^Ja_J$ and also condition 3(d) follows. \qed

\medskip
We also have

\begin{proposition}
\label{prop2.3}
Under
the assumptions of Proposition~\ref{prop1.3} there exists a unique graded
antilinear involution $*: \Gamma^{\w} \to \Gamma^{\w}$ such that

\begin{enumerate}
\begin{enumerate}
\item $(\th \wedge \th')^* = (-1)^{kl}\th'{}^* \wedge \th^*$, $\th \in
\Gamma^{\w k}$, $\th' \in \Gamma^{\w l}$,

\item $\Psi^{\w}(\th^*) = (\Psi^{\w}(\th))^{*\otimes *}$, $\th \in
\Gamma^{\w}$,

\item $d(\th^*) = d(\th)^*$, $\th \in \Gamma^{\w}$,

\item $*$ on $\Gamma^{\w 0}$, $\Gamma^{\w 1}$ coincides with the
original one.
\end{enumerate}
\end{enumerate}
\end{proposition}

\medskip
\noindent
{\bf Proof}.  We use the notation of the proof of Theorem~\ref{th2.2}.  We
define $*: \Gamma^{\otimes n} \to \Gamma^{\otimes n}$ by
\[
(\o_1 \otimes_{\cal C}
\ldots \otimes_{\cal C} \o_n)^* = (-1)^{\frac {n(n-1)}{2}}
\o_n^* \otimes_{\cal C} \ldots \otimes_{\cal C} \o_1^*
\]
for $\o_1,\dots,\o_n \in \Gamma^{\w 1}$.  Using $(4.14)$ of \cite{INH},
$(\o^{ij})^* = dx^j \otimes_{\cal C} dx^i + R^{ji}{}_{lk}dx^l \otimes_{\cal C}
dx^k = \o^{ji}$.  Hence $L^* \subset L$ and we get $*: \Gamma^{\w n} \to
\Gamma^{\w n}$ satisfying conditions (a)(d).  Condition (b)
follows from (\ref{eq2.2'}).  We see that conditions (a)(d) imply
\begin{equation}
\label{eq2.5}
(dx^Ia_I)^* = a_I^*(dx^I)^* = (-1)^{|I|(|I| - 1)/2}a_I^*dx^{I'},
\end{equation}
where $I' = (i_n,\dots,i_1)$ for $I = (i_1,\dots,i_n)$.  Therefore $*$ is
unique.  Set $\th = dx^Ia_I$.  According to (\ref{eq2.5}), (\ref{eq2.4}),
condition (a) and (\ref{eq2.3}),
\[
d(\th^*) = (-1)^{|I|(|I|-1)/2}d(a_I)^* \wedge dx^{I'} = ((-1)^{|I|}dx^I \wedge
da_I)^* = (d\th)^*
\]
and condition (c) follows. \qed

\medskip
In particular, quantum Minkowski spaces with ${\tilde F} = 0$ (described after
the proof of Theorem~\ref{th1.2}) admit a unique exterior algebra as described
in Theorem~\ref{th2.2} and Proposition~\ref{prop2.3}.  Moreover, using the
arguments of the proof of Theorem $1.9$ of \cite{POI}, we get

\begin{proposition}
\label{prop2.4}
Let $M$ be a quantum Minkowski space.  Then $\dim A_k = \four$,
$k = 0,1,2,3,4$, $\dim A_k = 0$ for $k > 4$.
\end{proposition}

\section{Differential operators}
\label{sec3}

Here we introduce and investigate the properties of the momenta
${P}_j = i\p_j$,
${P}^j = i\p^j$, the Laplacian $\Bo$ and the Dirac operator
${P}\!\!\!\!/$.  In particular, the momenta commute with
$\Bo$.  We heuristically show that
${P}^j,\Bo$ are hermitian.

Throughout the Section we deal with the exterior algebra on $M$ as described in
Theorems~\ref{th1.2} and \ref{th2.2}.  Further assumptions will be made later
on.  Let us recall that the partial derivatives $\p_i$ were defined by
(\ref{eq1.5}) and satisfy (\ref{eq1.8}), (\ref{eq1.6}).  Their values can be
computed using the unital homomorphism (\ref{eq1.5a}).  They satisfy the
following

\begin{proposition}
\label{prop3.1}
\begin{enumerate}
\item One has
\setcounter{equation}{0}
\begin{equation}
\label{eq3.1}
\p_l\p_k = R^{ij}{}_{kl} \p_j\p_i,\ i,j = 1,\dots,N.
\end{equation}

\item There exist $X_{i}{}^k$, $Y_i \in {\cal B}'$ such that
\begin{equation}
\label{eq3.7}
\r_{i}{}^k = (X_{i}{}^k \otimes id)\Psi,\ \p_i = (Y_i \otimes id)\Psi.
\end{equation}

\item One has
\begin{equation}
\label{eq3.7'}
\p_c\r_a{}^t=\r_d{}^t\p_bR^{bd}{}_{ac}, \quad R^{st}{}_{bd}\r_c{}^d\r_a{}^b=
\r_d{}^t\r_b{}^sR^{bd}{}_{ac}.
\end{equation}
\end{enumerate}
\end{proposition}

\medskip
\noindent
{\bf Proof}.  {\bf Ad 1}.
Let $a \in {\cal C}$.  Using (\ref{eq1.5}) and (\ref{eq2.2}), we
get
\[
\begin{array}{rll}
0 &= &dda = d(dx^i\p_i(a)) \\
&= &-dx^i \wedge d(\p_i(a)) \\
&= &-dx^i \wedge dx^j\p_j\p_i(a) \\
&= &-\frac {1}{2} (dx^i \wedge dx^j - R^{ij}{}_{kl}dx^k\wedge dx^l)
\p_j\p_i(a).
\end{array}
\]
But
\[
\frac {1}{2}(dx^i \wedge dx^j - R^{ij}{}_{kl}dx^k \wedge dx^l) =
(A_2)^{ij}{}_{kl} dx^k \wedge dx^l = (\a'_s)^{ij}\o^s
\]
(we have used (\ref{eq2.0'}), (\ref{eq2.1})).  In virtue of the
independence of
$\o^s$,
 $(\a'_s)^{ij}\p_j\p_i = 0$.  Multiplying by $(\alpha^s)_{kl}$ and using
(\ref{eq2.0'}), we get (\ref{eq3.1}).

\medskip
{\bf Ad 2}.  We set
\begin{equation}
\label{eq3.8}
X(a) = \left( \begin{array}{cc}
(f^l{}_{j}(S(a)))_{j,l=1}^N & 0 \\
0 & \e(a)
\end{array} \right),\ a \in {\cal A},
\end{equation}
\begin{equation}
\label{eq3.9}
X(y^i) = \left( \begin{array}{cc}
(Z^{il}{}_j)_{j,l=1}^N & (\d^{i}{}_j)_{j=1}^N \\
0 & 0
\end{array} \right),\ i = 1,\dots,N.
\end{equation}
There exists a unital homomorphism $X: {\cal B} \to M_{N+1}({\cb})$ which
satisfies (\ref{eq3.8})--(\ref{eq3.9})
($X(I_{\cal A}) = 1\!\!1$;
using $\tilde F=0$ and (2.18) of \cite{INH} for $b=v^k{}_l$, $v\in\mbox{ Rep }
H$, we show that $X(a)$,
$X(y^i)$ satisfy
(\ref{eqA})--(\ref{eqB})
for $a=v^i{}_j$, $v\in\mbox{ Rep }H$).  We have
\[
X = \left( \begin{array}{cc}
(X_{j}{}^l)_{j,l=1}^N & (Y_j)_{j=1}^N \\
0 & \e
\end{array} \right),
\]
where $X_{j}{}^l,Y_j \in {\cal B}'$.
By a direct computation (cf. the proof of
Theorem~\ref{th1.2})
\begin{equation}
\label{eq3.10}
{\cal L}(x) = (X \otimes id)\Psi(x)
\end{equation}
for $x = x^i \in {\cal C}$.  But ${\cal L}$, $X$ are both unital
homomorphisms.  Hence (\ref{eq3.10}) follows for all $x \in {\cal C}$. This
proves (\ref{eq3.7}).

\medskip
{\bf Ad 3.} We set $X_i{}^+=Y_i$, $X_+{}^i=0$, $X_+{}^+=\e$. We claim that
\begin{equation}
\label{eq3.10'}
(X_a{}^b * X_c{}^d)K_{bd}{}^{st} = K_{ac}{}^{bd}(X_b{}^s * X_d{}^t),\ a,c,s,t
= 1,\dots,N,+,
\end{equation}
where
\[
K = \left( \begin{array}{cccc}
R^T &0 &0 &0 \\
0 &0 &1\!\!1 &0 \\
0 &1\!\!1 &0 &0 \\
0 &0 &0 &1
\end{array} \right) .
\]
Indeed, since $X$
is a homomorphism, it is enough to
check (\ref{eq3.10'}) on $a \in {\cal A}$ (which follows from (\ref{eq3.8})
and the last formula before Proposition 3.14 of \cite{INH}),
and on $y_i$, $i = 1,\dots,N$, (which can be obtained  by a direct
computation -- see (3.61) of \cite{INH}). Using (\ref{eq3.10'}) and
(\ref{eq3.7}), one easily gets (\ref{eq3.7'}) and also (\ref{eq3.1}). \qed

\medskip
Let us notice that
\[
\begin{array}{rll}
dx^i\p_i(x^ka) &= &d(x^ka) = (dx^k)a + x^kdx^l\p_l(a) \\
&= &dx^i[{\d^k}_ia + R^{kl}{}_{in} x^n\p_l(a) + Z^{kl}{}_i \p_l(a)],
\mbox{ i.e.}
\end{array}
\]
\begin{equation}
\label{eq3.1''}
\p_ix^k = {\d^k}_i + (R^{kl}{}_{in} x^n + Z^{kl}{}_i)\p_l.
\end{equation}

In the following we assume that $g = (g^{ab})_{a,b=1}^N \in \mbox{
Mor}(I,\Lambda \otimes \Lambda)$ is a fixed invertible matrix.
It is called the metric tensor.  Then (see (1.12) of \cite{INH} and
(\ref{eq1.0'}))
$Rg =
g$.
Moreover $g^{-1} = (g_{ab})_{a,b=1}^N \in
\mbox{ Mor}(\Lambda \otimes \Lambda,I)$
($\Lambda g\Lambda^T=g$ implies $\Lambda^{-1}=g\Lambda^T g^{-1}$,
$\Lambda^T g^{-1}\Lambda=g^{-1}$).
  In the case with $*$ (as in
Propositions~\ref{prop1.3} and \ref{prop2.3})
one has ${\tilde g} \in \mbox{
Mor}(I, \Lambda \otimes \Lambda)$, where ${\tilde g}^{ij} = \overline{g^{ji}}$.
Then we also assume ${\tilde g} = g$.  Such $g$ exists e.g. for quantum
Poincar\'e groups and in this case is given (up to a nonzero real
multiplicative factor) by $q^{1/2}m$, where $m$ is given after Theorem~$1.12$
of \cite{POI} (we use (2.2) of \cite{POI},
$m^{-1}=(E'\otimes E'\tau)(\je\otimes X^{-1}\otimes\je)(V\otimes V)$).

We set $\p^a = g^{ab}\p_b$, ${P}_a = i\p_a$, ${P}^a = i\p^a$,
$\Bo =
g^{ij}\p_j\p_i = g_{ij}\p^i\p^j$.  Moreover, $\p^i\p^j = R^{ij}{}_{kl}
\p^k\p^l$
(we use (\ref{eq3.1}) and twice $g^{jb}R^{dc}{}_{ba}=R^{jd}{}_{ak}g^{kc}$,
which follows from $R^{jd}{}_{ak}=(R^{-1})^{jd}{}_{ak}=
f^d{}_a((\Lambda^{-1})^j{}_k)$ and $\Lambda^{-1}=g\Lambda^T g^{-1}$).
Using $(3.64)$ of \cite{INH} for $m = g$ and (\ref{eq3.1}), one easily
gets
\begin{equation}
\label{eq3.2}
\Bo \p_k = \p_k \Bo.
\end{equation}
Therefore $\Bo$ commutes with all $\p_k$, $\p^k$, ${P}_k$,
${P}^k$.  Moreover,
(\ref{eq1.6}) implies
\begin{equation}
\label{eq3.3}
(id \otimes \Bo)\Psi(a) = \Psi(\Bo(a)),
\end{equation}
\begin{equation}
\label{eq3.4}
(id \otimes \p^j)\Psi(a) = ((\Lambda^{-1})^j{}_i \otimes I)\Psi(\p^i(a)).
\end{equation}
The momenta ${P}_j$, ${P}^j$ transform under the action of the
inhomogeneous
quantum group $G$
in the same way
as the partial derivatives $\p_j$, $\p^j$ in (\ref{eq1.6}),
(\ref{eq3.4}).

The Dirac gamma matrices are defined as matrices
$\g^a\in M_d(\cb)$, $a =
1,\dots,N$, satisfying
\begin{equation}
\label{eq3.5}
\g^a\g^b + R^{ba}{}_{dc}\g^c\g^d = 2g^{ba}1\!\!1,\ a,b = 1,\dots,N,
\end{equation}
(cf. \cite{RTF}, \cite{Sch}).
It seems
that one can analyze irreducible representations of
the relations (\ref{eq3.5})
only case by case for particular matrices $R$.  In the following we assume that
such an irreducible representation has been chosen.  Then the Dirac operator
$\p\!\!\!/ = \p_a\g^a$ satisfies
\[
\begin{array}{rll}
{\p\!\!\!/}^2 &= &\p_c\p_d\g^c\g^d = \frac {1}{2}(\p_c\p_d + R^{ba}{}_{dc}
\p_a\p_b)
\g^c\g^d \\
&= &\frac {1}{2} \p_a\p_b(\g^a\g^b + R^{ba}{}_{dc}\g^c\g^d) \\
&= &g^{ba}\p_a\p_b \cdot 1\!\!1 = \Bo \cdot 1\!\!1
\end{array}
\]
(cf. \cite{Sch}).  Let ${P}\!\!\!\!/ = i\p\!\!\!/$. Then
${{P}\!\!\!\!/}^{\ 2}=-\Bo\cdot\je$.

Later on we shall need the following

\begin{proposition}
\label{prop3.2}
One has
\begin{equation}
\label{eq3.6}
\p^jx^k = g^{jk} + R^{jk}{}_{ab} x^a\p^b - (RZ)^{jk}{}_b \p^b
\end{equation}
\end{proposition}

\medskip
\noindent
{\bf Proof}. According to (\ref{eq3.1''}), $\p^jx^k = g^{jk} +
G^{jk}{}_{ab}x^a\p^b + V^{jk}{}_b\p^b$, where
\[
\begin{array}{rll}
G^{jk}{}_{ab} &= &g^{ji}R^{kl}{}_{ia}g_{lb} \\
&= &f^k{}_a(g^{ji}\Lambda^l{}_ig_{lb}) \\
&= &f^k{}_a((\Lambda^{-1})^j{}_b) \\
&= &(R^{-1})^{jk}{}_{ab} \\
&= &R^{jk}{}_{ab}, \\
V^{jk}{}_b &= &g^{ji}Z^{kl}{}_i g_{lb} \\
&= &\eta^k(g^{ji} \Lambda^l{}_ig_{lb}) \\
&= &\eta^k((\Lambda^{-1})^j{}_b) \\
&= &-(RZ)^{jk}{}_b
\end{array}
\]
(cf. the proof of Proposition $4.5.2$ of \cite{INH}).\qed

\medskip
In the remaining part of the section we shall heuristically prove that the
operators ${P}^a,\Bo$ are hermitian.  A strict
proof
of that fact would need the existence of topological structures on the quantum
spaces
$G,M$.  Therefore the considerations below serve as a motivation for such a
topological approach.  We shall assume the existence of a $G$-invariant measure
$\mu$ defined on some ``functions'' $x$ on $M$ (elements of ${\cal C}$ are also
``functions'' on $M$ but we don't expect $\mu$ to act on them; nevertheless we
shall disregard such subtleties here).  Thus $(id \otimes
\mu)\Psi x = \mu(x)I_{\cal B}$.  One has the scalar product $(x \mid y) =
\mu(x^*y)$.  Using (\ref{eq3.7}), we obtain
\[
\begin{array}{rll}
(I \mid \p_ix) &= &\mu(\p_ix) \\
&= &\mu(Y_i \otimes id)\Psi(x) \\
&= &Y_i[(id\otimes \mu)\Psi(x)] \\
&= &Y_i(I_{\cal B})\mu(x) \\
&= &0\
\end{array}
\]
$(0 = \p_iI = Y_i(I)I)$,
$(I \mid {P}^ax) = 0$, $({P}^a)^*I = 0$.
Also $P^aI=0$. Moreover, hermitian
conjugate of
(\ref{eq3.6}) yields
\[
x^k(\p^j)^* = g^{kj} + R^{kj}{}_{ba}(\p^b)^*x^a + Z^{kj}{}_b (\p^b)^*,
\]
i.e.
\[
(\p^b)^*x^a = -g^{ba} + R^{ba}{}_{kj} x^k(\p^j)^* - (RZ)^{ba}{}_s(\p^s)^*.
\]
This means that ${P}^b$ and $({P}^b)^*$ commute with $x^a$ in the
same way.  Thus
they act in the same way on all elements of ${\cal C}$ and ${P}^b$ is
hermitian.
This and $g_{ab} = \overline{g_{ba}}$ show that
$\Bo$ is also hermitian.

\section{Solutions of Klein--Gordon and Dirac equations}
\label{sec4}

We shall consider formal solutions $\var$ of Klein--Gordon equation $(\Bo +
m^2)\var = 0$ and Dirac equation ${P}\!\!\!\!/ \var = m\var$ obtained
using the plane waves $e^{-ip_ax^a}$.
But now $p_a$ are not (in a general case)
eigenvalues of ${P}_a$.

Setting $F_i{}^l(x^k) = R^{kl}{}_{in}x^n + Z^{kl}{}_i$ and using $\p_ix^k =
{\d^k}_i + F_i{}^l(x^k)\p_l$ (see (\ref{eq3.1''})), one gets
\setcounter{equation}{0}
\begin{equation}
\label{eq4.1}
\p_j(x^{k_1}\dots x^{k_n}) = \sum_{m=0}^{n-1}
F_j{}^{i_1}(x^{k_1})F_{i_1}{}^{i_2}(x^{k_2})\dots
F_{i_{m-1}}{}^{k_{m+1}}(x^{k_m})x^{k_{m+2}}\dots x^{k_n},
\end{equation}
where $m = 0$ term reads ${\d^{k_1}}_jx^{k_2}\dots x^{k_n}$.
We shall consider two cases:
$Z = 0$ and $R = \tau$.

1.\ \  $Z = 0$.  Let ${\cal F}$ be a new ${}^*$-algebra generated by $p^k$,
$k =
1,\dots,N$, satisfying $(p^k)^* = p^k$ and $p^kp^l = R^{lk}{}_{ji}p^ip^j$ (same
relations as for ${P}^k$ but w.r.t. the opposite multiplication).
Set $p_a = g_{ab}p^b$, $\pe=p_j\otimes\g^j$.
Then
\begin{equation}
\label{eq4.2}
p_ap_b = p_kp_l R^{kl}{}_{ab}.
\end{equation}
We put $x \otimes p = x^a \otimes p_a \in {\cal C} \otimes {\cal F}$.
So in a sense we assume that $p_a$ commute with all $x^b$ and $\p_m$.
This doesn't seem to be the case in \cite{MajidE}, \cite{Ch}, \cite{Meyer},
\cite{Kulish} (cf. also \cite{CSW}, \cite{Pillin}, \cite{ChD}) although
calculations in the present case 1. (up to the moment when
we use $\pi$) are quite similar as in these papers.
Using
(\ref{eq4.1}) and (\ref{eq4.2}), one obtains
\begin{equation}
\label{eq4.3}
(\p_j \otimes id)(x \otimes p)^n = \sum_{m=0}^{n-1}
(I\otimes p_j)(x\otimes p)^{n-1} = n(I
\otimes p_j)(x \otimes p)^{n-1}.
\end{equation}
Thus in the sense of formal power series w.r.t. $t$
\begin{equation}
\label{eq4.3'}
(\p_j \otimes
id)e^{-it(x\otimes p)} = -it(I \otimes p_j)e^{-it(x\otimes p)},
\end{equation}
\begin{equation}
\label{eq4.4}
(\Bo \otimes id)e^{-it(x\otimes p)} = -t^2(I \otimes s)e^{-it(x \otimes p)},
\end{equation}
where $s = g^{ij}p_ip_j = g_{im}p^mp^i = p_ip^i$ is a central, hermitian
element of ${\cal F}$.  Let $\pi$ be an irreducible ${}^*$-representation of
${\cal F}$ in a Hilbert space $H$ with an orthonormal basis $e_k$, $k \in K$.
We denote $\pi_{kl}(x) = (e_k \mid \pi(x)e_l)$, $x \in {\cal F}$.
In the remaining part of $Z=0$ case we assume that all performed
operations are well defined and have ``good'' properties. Acting $id
\otimes \pi_{kl}$ on (\ref{eq4.4}) and setting $\pi(s) = m^2 \in
{\rb}$, $(id \otimes \pi_{kl})(e^{-it(x \otimes p)}) = \var_{kl}^{(t)}$, one
obtains $\Bo \var_{kl}^{(t)} = -t^2m^2\var_{kl}^{(t)}$.  In a topological
approach one should be able to put $t = 1$. Then $\var =
\var_{kl}^{(t=1)}$ is a solution
of Klein--Gordon equation $\Bo \var = -m^2{\varphi}$
in a usual sense (if $m \ge 0$).  Also (see (\ref{eq4.3'}))
\[
{P}^j\var_{kl}^{(t=1)} = \pi_{ks}(p^j)\var_{sl}^{(t=1)}
\]
so one should get a real spectrum of $P^j$ ($p^j$ are selfadjoint).

Let us pass to the Dirac equation.  We denote the canonical basis in ${\cb}^d$
 by $\e_m$, $m = 1,\dots,d$, $\var_{klm}^{(t)} = \var_{kl}^{(t)}
\otimes \e_m$.  Tensoring (\ref{eq4.3}) by $\g^j\e_m$, one obtains
\[
{\p\!\!\!/}_{(13)}[(x \otimes p)^n \otimes \e_m] = n {\pe}_{(23)}[(x
\otimes p)^{n-1} \otimes \e_m],
\]
where e.g. $(13)$ means that $\p\!\!\!/$ acts in the first and the third
factors of the tensor product.  Thus
\[
{{P}\!\!\!\!/}_{(13)}[e^{-it(x \otimes p)} \otimes \e_m] =
t{\pe}_{(23)}[e^{-it(x \otimes p)} \otimes \e_m].
\]
Acting by $id \otimes \pi_{kl} \otimes id$, we get
\begin{equation}
\label{eq4.5}
{P}\!\!\!\!/ {\varphi}_{klm}^{(t)} = t\pi_{ks}(p_a)(\g^a)^i{}_m
{\varphi}_{sli}^{(t)}.
\end{equation}
Set $U^{si}{}_{kn} = \pi_{ks}(p_a)(\g^a)^i{}_n$.  One has $U^2 = m^21\!\!1$.
The possible eigenvalues of $U$ are $\pm m$.  Suppose $Uv = mv$ for $v
\in H \otimes {\cb}^d$.  Put ${\varphi}_{vl}^{(t)} =
v^{si}\var_{sli}^{(t)}$.  Then ${P}\!\!\!\!/ {\varphi}_{vl}^{(t)} =
tm{\varphi}_{vl}^{(t)}$.  In a topological approach one should be able to
put $t = 1$ and find a solution $\var = \var_{vl}^{(t=1)}$ of Dirac equation
${P}\!\!\!\!/ \var = m\var$ (if $m \ge 0$).

2. \ \  $R = \tau$, i.e. $R^{ij}{}_{kl} = {\d^i}_l {\d^j}_k$.
Then $Rg = g$ implies
$g^{ab} = g^{ba} = \overline{g^{ab}}$.  According to the proof of Proposition
$4.5.2$ of \cite{INH}, $RZ = -\tau{\bar Z}$, $Z = -{\bar Z}$, $Z^{ab}{}_c \in
i{\rb}$.  Let $p^k \in {\rb}$, $k = 1,\dots,N$, $p_j = g_{jk}p^k \in
{\rb}$.  Setting ${U}_i{}^l = Z^{kl}{}_i p_k \in i{\rb}$, $x \cdot p
= x^jp_j \in {\cal C}$, we get $F_i{}^l(x \cdot p) =
{\d^l}_i(x \cdot p) + {
U}_i{}^l$.  This and (\ref{eq4.1}) yield $(({U}^k)_i{}^j = {
U}_i{}^{i_1} \cdot {U}_{i_1}{}^{i_2} \cdot \dots \cdot {
U}_{i_{k-1}}{}^j)$
\[
\p_j(x \cdot p)^n =
\sum_{ \begin{array}{c} a_0 + \ldots + a_r = m - r \\ r \ge 0,\ 0 \le m \le
n-1\end{array}}
[{U}^{a_0}(x \cdot p){U}^{a_1}(x \cdot p) \cdot \dots
\cdot (x \cdot p){U}^{a_r}]_j{}^b p_b(x \cdot p)^{n-m-1}.
\]
Let $G = \sum_{s=0}^{\infty} (tU)^s$
(in the sense of formal power series w.r.t. t).
Since $x \cdot p$ commutes with ${
U}$
\[
\p_j \sum_{n=0}^{\infty} (tx \cdot p)^n = \sum_{l=0}^{\infty} \sum_{r=0}^l
(G^{r+1})_j{}^b tp_b(tx \cdot p)^l.
\]
But using the mathematical induction
\[
\sum_{r=0}^l G^{r+1} = \sum_{k=1}^{\infty}
{(\hskip-5pt\begin{array}{c}{\scriptstyle l+k}\\[-8pt]%
{\scriptstyle k}\end{array}\hskip-5pt)}
(tU)^{k-1}.
\]
Therefore
\[
\p_j(x \cdot p)^n = \sum_{k=1}^n
{(\hskip-5pt\begin{array}{c}{\scriptstyle n}\\[-8pt]%
{\scriptstyle k}\end{array}\hskip-5pt)}
(U^{k-1})_j{}^b p_b(x \cdot p)^{n-k},
\]
\begin{equation}
\label{eq4.6}
\begin{array}{rll}
\p_j e^{-itx\cdot p} &= &[\rho(-itU)]_j{}^b (-itp_b)e^{-itx\cdot p},
\end{array}
\end{equation}
where $\rho(x) = \sum_{k=1}^{\infty} \frac {x^{k-1}}{k!} = \frac {e^x-1}{x}$.
Using $(3.62)$ of \cite{INH} for $m = g$, one obtains
\begin{equation}
\label{eq4.6'}
Z^{kl}{}_r g^{rj} = -Z^{kj}{}_s g^{ls},
\end{equation}
${U}_r{}^lg^{rj} = -g^{ls}U_s{}^j$.  In virtue of the mathematical
induction $({U}^m)_r{}^l g^{rj} = g^{ls}[(-U)^m]_s{}^j$, $\rho(-itU)_r{}^l
g^{rj} = g^{ls}\rho(itU)_s{}^j$.  Thus
\[
\Bo e^{-itx\cdot p} = g^{mj} \rho(-itU)_m{}^a \rho(-itU)_j{}^b(-t^2p_ap_b)
e^{-itx\cdot p} = g^{as}h(-itU)_s{}^b(-t^2p_ap_b)e^{-itx\cdot p},
\]
 where $h(x)
= \rho(-x)\rho(x) = \left( \frac {\sinh(x/2)}{x/2} \right)^2$.
Set
\begin{equation}
\label{eq4.7}
m^2 = g^{as}h(-iU)_s{}^bp_ap_b \in {\rb}.
\end{equation}
In a
topological approach one should be able to put $t = 1$ and find a solution
$\var = e^{-ix\cdot p}$ of Klein--Gordon equation $(\Bo + m^2)\var = 0$
(if $m \ge 0$).

Using (\ref{eq4.6}), one would also have
\[
P^j{\varphi} = {\cal P}^{j(t=1)}{\varphi},
\]
where ${\cal P}^{j(t=1)} = g^{js}\rho(-iU)_s{}^b p_b \in {\rb}$.

Let us pass to the Dirac equation.  Setting $\var_r = e^{-itx\cdot p} \otimes
\e_r$, one gets ${P}\!\!\!\!/ \var_r = ({\cal P}\!\!\!\!/)^k{}_{r} \var_k$,
where ${\cal P}\!\!\!\!/ = {\cal P}_j\g^j$, ${\cal P}_j =
\rho(-itU)_j{}^b tp_b \in {\rb[[t]]}$.
In a
topological approach one should be able to put $t = 1$.
 Then
${{\cal P}\!\!\!\!/}^2 = g^{js}{\cal P}_j{\cal P}_s1\!\!1$,
the possible eigenvalues of ${\cal
P}\!\!\!\!/$ are $\pm m$, $m = \sqrt{g^{js}{\cal P}_j{\cal
P}_s}$.  Suppose ${\cal P}\!\!\!\!/ v = mv$ for $v \in {\cb}^d$.
Put $\var = v^r\var_r=e^{-ix\cdot p}\otimes v$.
Then
$\var$ is a solution of Dirac equation ${P}\!\!\!\!/ \var = m\var$
(if $m \ge 0$).

Let us notice that spectra of ${P}^j$, $\Bo$ obtained
(formally) above are real in agreement with hermiticity of those operators,
heuristically asserted at the end of Section~\ref{sec3}.

The case $R = \tau$ covers e.g. the
quantum Poincar\'e groups of case 1), $s = 1$,
$t = 1$, $t_0 = 0$. Many of them are described by Remark $1.8$ of \cite{POI}
and
\cite{Z}.  Then eigenvalues $\lambda$ of $-iU$ are real or imaginary which
yields factors $\left( \frac {\sinh(|\lambda|/2)}{|\lambda|/2} \right)^2$,
$\left(
\frac {\sin(|\lambda|/2)}{|\lambda|/2} \right)^2$ in (\ref{eq4.7}) (when $U$ is
diagonalizable).
Despite the sine factors,
 $m^2$ of (\ref{eq4.7}) is
never bounded in the cases given by \cite{Z}.
Hyperbolic sine factors can improve the properties of the propagator
$\frac {i}{g^{as}h(-iU)_s{}^bp_ap_b-M^2}$ ($M \ge 0$)
(although it is never integrable in the cases given by \cite{Z})
e.g. $b=\e e_1\w H$, $a=0$, $c=0$ in \cite{Z} ($\e\in\rb$) corresponds to
$U_{03}=U_{30}=i\e p_1/2$, other $U_{ab}$ are $0$ and we get the propagator
\[
\frac{i}{(p_0^2-p_3^2)\left[\frac{\mbox{ sinh$(\e p_1/4)$}}
{\mbox{$\e p_1/4$}}
\right]^2-p_1^2-p_2^2-M^2}.
\]
A detailed analysis of Feynman diagrams is not yet possible since the measure
in the $p$--space and interaction factors are not known.
Let us also notice that for $H=SO(4)$,
 $U$ is antisymmetric and we can have only
factors of the second type.

\medskip
{\bf Remark.} In \cite{Sitarz} differential calculi corresponding to
a quantum Minkowski space of \cite{Lukierski} are considered. Despite
other choices of axioms the result is the same as here: there are no
$4$--dimensional covariant differential calculi in that case (due to
\cite{Z} $t_0\neq 0$ for \cite{Lukierski}).
 But nevertheless a different
approach suggests a similar form for the propagator as above \cite{Lukierski}.

\section{The Fock space}
\label{sec5}

Here we define the Fock space for noninteracting particles on $M$.  We assume
that $G$ has $CT$ Hopf ${}^*$-algebra structure
${\cal R}$ as in Theorem $3.1.3$--$4$ of
\cite{CQT} which is the case e.g. for quantum Poincar\'e groups of cases 1) 2)
3) 4) (except 1), $s = 1$, $t = 1$, $t_0 \ne 0$ and 4), $s = 1$, $b \ne 0$) as
proved in Theorem 3.2.3 of \cite{CQT}.

We follow the scheme of \cite{QMQS} but now we have the left (instead of right)
action.  The particles interchange operator $K: {\cal C} \otimes {\cal C} \to
{\cal C} \otimes {\cal C}$ is defined by
\[
K(x \otimes y) = {\cal
R}(y^{(1)} \otimes x^{(1)})(y^{(2)} \otimes x^{(2)})
\]
where
$\Psi(x) = x^{(1)} \otimes x^{(2)}$, $\Psi(y) =
y^{(1)} \otimes y^{(2)}$.  We set $K^{(m)} = {1\!\!1}^{\otimes (m-1)} \otimes K
\otimes {1\!\!1}^{\otimes (n-m-1)}: {\cal C}^{\otimes n} \to {\cal C}^{\otimes
n}$, $m = 1,2,\dots,n-1$.  Then one obtains the representation of the
permutation group
\[
\pi: \Pi_n \ni t_m = (m,m+1) \to \pi_{(m,m+1)} = K^{(m)}
\]
acting in ${\cal C}^{\otimes n}$ which defines the boson subspace ${\cal
C}^{\otimes_s n}$.  The group $G$ acts on ${\cal C}^{\otimes n}$ by the linear
mapping $\Psi^{\otimes n}: {\cal C}^{\otimes n} \to {\cal B} \otimes {\cal
C}^{\otimes n}$ defined by
\[
\Psi^{\otimes n}(x_1 \otimes \ldots \otimes x_n) = x_1^{(1)} \cdot \ldots \cdot
x_n^{(1)} \otimes x_1^{(2)} \otimes \ldots \otimes x_n^{(2)},
\]
$x_1,\dots,x_n \in {\cal C}$.  One has $\Psi^{\otimes n} K^{(m)} = (id \otimes
K^{(m)})\Psi^{\otimes n}$, i.e. the actions of $G$ and $\Pi_n$ agree.  In
particular, one can restrict $\Psi^{\otimes n}$ to ${\cal C}^{\otimes_s n}$
getting $\Psi^{\otimes_s n}: {\cal C}^{\otimes_s n} \to {\cal B} \otimes {\cal
C}^{\otimes_s n}$.  If $W: {\cal C} \to {\cal C}$ is an operator related to a
single particle then the corresponding $n$-particle operator is given by
\[
W^{(n)} = \sum_{m=1}^n \pi_{(1,m)} ( W\otimes {1\!\!1}^{\otimes
(n-1)})\pi_{(1,m)} =
\frac{1}{(n-1)!}\sum_{\sigma\in\Pi_n} \pi_{\sigma} (W\otimes {1\!\!1}^{\otimes
(n-1)})\pi_{\sigma}^{-1}:
{\cal C}^{\otimes_s n} \to {\cal C}^{\otimes_s n}
\]
(the $m$-th term in the first sum is the operator in ${\cal C}^{\otimes n}$
corresponding to the $m$-th particle).
We can also define the Fock space $F = \oplus_{n=0}^{\infty} {\cal
C}^{\otimes_s n}$ and the operator $\oplus_{n=0}^{\infty} W^{(n)}$ acting in
$F$.

For particles of mass $m$ we should consider $\ker(\Bo + m^2)$ instead of
${\cal C}$ and a scalar product there but it would lead us beyond the scope
of the present article (heuristically e.g. $W=P^k$, $k=1,\ldots,N$,
would be hermitian operators in such a space -- see (\ref{eq3.2})).

\bigskip
\begin{center}
{\bf Acknowledgments}
\end{center}

\medskip
I am grateful to Professor W. Arveson and other faculty members for
their kind hospitality at UC Berkeley.  The author is thankful to Professors J.
Lukierski, N. Reshetikhin, B. Zumino, Dr W. Biedrzycki and Dr A. Schirrmacher
for fruitful discussions.

\end{document}